\documentclass{article}
\usepackage[dvips]{graphicx}
\usepackage{amsthm,amssymb,amsfonts,epsfig}

\newtheorem{theorem}{Theorem}

\newtheorem{corollary}[theorem]{Corollary}

\newtheorem{lemma}[theorem]{Lemma}

\begin{document}

\title{On the heterochromatic number of hypergraphs associated to geometric graphs and to matroids\thanks{Partially supported by Conacyt, M{\'e}xico.} } 

\author{Juan Jos{\'e} Montellano-Ballesteros\\ 
Instituto de Matem{\'a}ticas,\\
Universidad Nacional Aut{\'o}noma de M{\'e}xico.\\
Area de la Investigaci{\'o}n Cient{\'i}fica,\\
Circuito Exterior, Ciudad Universitaria,\\
Coyoac{\'a}n,  D. F., 04510,\\
M{\'e}xico.\\
\texttt{juancho@matem.unam.mx}
\and
Eduardo Rivera-Campo\\
Departamento de Matem{\'a}ticas,\\
Universidad Aut{\'o}noma Metropolitana-Iztapalpa.\\
Av. San Rafael Atlixco 186, Col. Vicentina,\\
Iztapalapa, D.F., 09340,\\
M{\'e}xico.\\
\texttt{erc@xanum.uam.mx}}

\date{}

\maketitle

\begin{abstract}

The heterochromatic number $h_{c}\left( H\right)$ of a non-empty hypergraph $H$ is the smallest integer $k$ such that for every colouring of the vertices of $H$ with exactly $k$ colours, there is a hyperedge of $H$ all of whose vertices have different colours. We denote by  $\nu \left( H\right)$ the number of vertices of $H$ and by  $\tau \left( H\right)$ the size of the smallest set containing at least two vertices of each hyperedge of $H$.  For a complete geometric graph $G$ with $n\ge3$ vertices let  $H=H \left( G\right)$ be the hypergraph whose vertices are the edges of $G$ and whose hyperedges are the edge sets of plane spanning trees of $G$. We prove that if $G$ has at most one interior vertex, then $h_{c} \left(H\right)= \nu \left( H\right) - \tau \left( H\right) + 2$. We also show that $h_{c} \left( H\right)  = \nu \left( H\right) - \tau \left( H\right) + 2$ whenever $H$ is a hypergraph with vertex set and hyperedge set given by the ground set and the bases of a matroid, respectively.

\end{abstract}


\section{Introduction}
 
Let $k$ be a positive integer and $H$ be a hypergraph with at least $k$ vertices. A \emph{k-colouring} of  $H$ is an assigment of colours to the vertices of $H$ that uses exactly $k$ colours.  Given a colouring $c$ of a hypergraph $H$, a \emph{heterochromatic hyperedge} is a hyperedge $e$ of $H$ such that $c$ assigns different colours to different vertices of $e$.

The \emph{heterochromatic number} $h_{c}\left( H\right)$ of  a non-empty hypergraph is the smallest integer $k$ such that $H$ contains a heterochromatic hyperedge for each $k$-colouring of $H$. The heterochromatic number was defined by Arocha \emph {et   al}  \cite{ABN}  and is closely related to \emph{anti-Ramsey numbers}  \cite{ESS}  and to the \emph{upper chromatic number} of mixed hypergraphs  \cite{V}.

A \emph{double transversal of hyperedges} of a loopless hypergraph $H$ is a set $T$ of vertices of $H$ such that each hyperedge of $H$ contains at least two vertices in $T$. We denote by  $\nu \left( H\right)$ and $\tau \left( H\right)$, the number of vertices and the size of the smallest double transversal of hyperedges of $H$, respectively.

A general lower bound for the heterochromatic number of loopless  non-empty hypergraphs is obtained as follows: Consider a double transversal $T$ of hyperedges of $H$ with $\left\vert T \right\vert = \tau \left( H\right)$. Assign  colour 1 to every vertex in $T$ and a different colour to each of the remaining $\nu \left( H\right) - \tau \left( H\right)$ vertices of $H$.  Since $T$ is a double transversal of hyperedges of $H$, there are no heterochromatic edges of $H$ for this $\left( \nu \left( H\right) - \tau \left( H\right) + 1 \right)$-colouring of $H$ and therefore $h_{c} \left( H\right)   \geq \nu \left( H\right) - \tau \left( H\right) + 2$.

Jiang and West \cite{JW} proved that $h_{c} \left( H_{n}\right)={n-2 \choose 2} + 2$ if $H_{n}$ is the hypergraph whose vertices are the edges of a complete graph $G$ with $n$ vertices and whose hyperedges are the the edge sets of all spanning trees of $G$. Notice that, in this case, $\nu \left( H_{n}\right)={n\choose 2}$, $ \tau \left( H_{n}\right)=2n-3$ and $h_{c} \left( H_{n}\right)  = \nu \left( H_{n}\right) - \tau \left( H_{n}\right) + 2$.

In this article we study  hypergraphs $H$ associated to complete geometric graphs and to matroids for which the general lower bound also give the exact value of $h_{c} \left( H\right)$. For a complete geometric graph $G$, we associate a hypergraph $H\left( G \right)$ with vertex set given by the edges of $G$ and whose hyperedges are the sets of edges of non-self intersecting spanning trees of $G$. For a matroid $M$ we consider the hypergraph  $H\left( M \right)$  whose vertices are the elements of the ground set of $M$ and whose hyperedges are the bases of $M$. 

We show that $h_{c} \left( H\right)  = \nu \left( H\right) - \tau \left( H\right) + 2$ if $H$ is the hypergraph associated as above to  a complete  geometric graph with at most one interior vertex or to a matroid.

\section{Geometric graphs}

Let $P$ be a set of points in general position in the plane, a  \emph{geometric graph} on $P$ is a graph $G$ with vertex set $P$ drawn in such away that each edge is a straight line segment with ends in $P$.  

The \emph{complement}  of a geometric graph $G$ is the geometric graph $G^{c}$ with the same vertex set $P$  whose edges are all the line segments with ends in $P$ which are not edges of $G$.

A \emph{plane spanning tree} of a geometric graph $G$ is a non-selfintersecting subtree that contains every vertex of $G$.

A tree $T$ is a \emph{caterpillar} if the removal of the endpoints of $T$ leaves a path called the  \emph{body} of $T$. Let $P$ be a set of points in general position in the plane. A plane geometric tree $R$ with vertex set $P$ is a \emph{geometric caterpillar} if $R$ is a caterpillar such that the entire body of $R$ lies in the boundary of  the convex hull $CH\left( P\right)$ of $P$ and  for each leg $e$ of $R$, the straight line containing $e$ does not intersect $R$ at any point not in $e$.

The following results, due to K\'{a}rolyi \emph{et al} \cite{ KPT} and to Urrutia-Galicia \cite{U}, respectively, give sufficient conditions for the complement of a geometric graph to contain a plane spanning tree.

\begin{theorem}
\label{karolyi}

If the edge set of a complete geometic graph $G$ is partitioned into two sets, then there exists a plane spanning tree of $G$, all of whose edges lie in the same part.

\end{theorem}

\begin{theorem}
\label{urrutia}

Let $R$ be a plane spanning tree of a complete geometric graph $G$. The complementary geometric graph $R^c$ contains a plane spanning tree if and only if $R$ is neither a star or a geometric caterpillar.

\end{theorem}

For any set $P$ of points in general position in the plane, we denote by  $i \left( P\right)$ the number of points of $P$ not lying in the boundary of $CH\left( P\right)$. 

\bigskip 

Lemma  \ref{transversal1} (due to  Garcia  \emph{et al} \cite{GHHNT})  and Lemma \ref{transversal2}  provide double transversals of plane spanning trees in complete geometric graphs $G$  with  $i \left( V \left( G\right) \right) \leq 1$. 

\begin{lemma}
\label{transversal1}

Let $P$ be a set of  $n\geq3$  points in convex position in the plane. If $R$ is a plane spanning tree of the complete geometric graph with vertex set $P$, then at least two edges of $R$ lie in the boundary of $CH\left( P\right)$. 

\end{lemma}

\begin{lemma}
\label{transversal2}

Let $P$ be a set of $n$ points in general position in the plane such that $i \left( P\right) = 1$ and $w$ be the unique point in $P$ not lying in  the boundary of $CH\left( P\right)$. Let  $Q$ be the set of edges of the boundary of $CH\left( P\right)$ together with two edges $uw$ and $vw$ such that the angle $\angle uwv$ is maximal. If $R$ is a plane spanning tree of the complete geometric graph with vertex set $P$, then at least two edges of $R$ lie in $Q$.

\end{lemma} 

\begin{proof}

Let $G$ be the complete geometric graph with vertex set $P$. If $n = 4$, then each spanning tree of $G$ has 3 edges and $Q$ contains all but one of the edges of $G$.  We proceed by induction assuming $n\geq5$ and that the result holds for every subset $P^\prime$ of $P$ with $i \left( P^\prime \right) = 1$. 

Let $R$ be a plane spanning tree of $G$. Either  $R$ contains two edges $ux$ and $vy$ in $Q$ which are incident in $u$ and $v$, respectively, or $R$ contains an edge $e$ incident in $u$ or in $v$ which is a diagonal of the boundary of  $CH\left( P\right)$. 

Supose  $R$ contains such a diagonal edge $e$. As  $\angle uwv$ is maximal, then $e$ cannot intersect the edges $uw$ and $wv$ in a point other than $u$ or $v$.  Let  $P^{-}$ be the set of points in $P$ lying on or to the left of $e$ and  $P^{+}$ be the set of points in $P$ lying on or to the right of $e$, notice that $w$ is an interior point of $P^{-}$ or an interior point of $P^{+}$. Again without loss of generality, assume $i \left( P^{-}\right) = 0$ and $i \left( P^{+}\right) = 1$.

Let $Q^{-}$ be the set of edges in the boundary of  $CH\left( P^{-}\right)$ and  $Q^{+}$ be the set of edges of the boundary of $CH\left( P^{+}\right)$ together with the two edges $uw$ and $vw$. By Lemma  \ref{transversal1}, the subtree $R^{-}$ of $R$ with vertex set $P^{-}$ contains two edges in $Q^{-}$ and by inducction, the subtree $R^{+}$ of $R$ with vertex set $P^{+}$ contains two edges in $Q^{+}$. Since $Q^{-} \cup Q^{+} = Q$ and  $Q^{-} \cap Q^{+} =\left\{ e\right\}$, at least two edges of $R$ lie in $Q$.

\end{proof}

We can now prove the main results of this section.

\begin{theorem}
\label{geometricas1}

Let  $G$ be a complete geometric graph with $n\geq3$ vertices. If $c$ is a colouring of the edges of $G$ with exactly  ${n \choose 2}  - n + 2$ colours, then $G$ has a heterochromatic plane spanning tree.

\end{theorem} 

\begin{proof}

Let $X$ be a heterochromatic set with ${n \choose 2}  - n + 2$ edges of $G$ and let $Y = E\left( G\right) \setminus X$. Since $\left\vert Y\right\vert  = n - 2$, no spanning tree of $G$ has all edges in $Y$. By Theorem \ref{karolyi}, graph $G$ has a plane spanning tree $R$ all of whose edges lie in $X = E\left( G\right) \setminus Y$. Clearly $R$ is a heterochromatic tree.

\end{proof}

\begin{theorem}
\label{geometricas2}

Let  $G$ be a complete geometric graph with $n\geq3$ vertices.  If $i  \left( V \left( G \right) \right) = 1$ and $c$ is a colouring of the edges of $G$ with exactly  ${n \choose 2}  - n + 1$ colours, then $G$ has a heterochromatic plane spanning tree.

\end{theorem} 

\begin{proof}

Let $X$ be a heterochromatic set with ${n \choose 2}  - n + 1$ edges of $G$ and let $Y = E\left( G\right) \setminus X$.  Since $\left\vert Y\right\vert  = n - 1$, either $Y$ is the set of edges of a plane spanning tree of $G$ or no spanning tree of $G$ has all edges in $Y$. As in the proof of Theorem \ref{geometricas1}, in the later case, $G$ has a heterochromatic plane spanning tree $R$.

Assume $Y$ is the set of edges of a plane spanning tree $S$ of $G$.  By Theorem \ref{urrutia}, either $S$ is a geometric caterpillar, $S$ is a star or $S^c$ contains a plane spanning tree $R$. In the latter case, $R$ is a heterochromatic plane spanning tree of $G$. 

For the case where $S$ is a geometric caterpillar but not a star let $y=uv$ be an edge in the body of $S$ with $d_{S} \left( u\right) \geq 2$ and $d_{S}\left( v\right) \geq 2$. There is an edge $x\in X$ with $c \left( x \right) =  c \left( y \right)$ since each colour is assigned to an edge in $X$.  

Let $S^\prime = \left( S  -  y  \right) + x $ and notice that  $S^\prime $ is not a star  by the choice of $y$. Suppose  $S^\prime $ is a geometric caterpillar in which case $x$ must lie in the body of $S^\prime $ since $y$ lies in the body of $S$. This implies that both $S$ and $S^\prime $ are paths whose union Á is the boundary of $CH\left( P \right)$ which is not possible since $i  \left( V \left( G \right) \right) = 1$. Therefore $S^\prime$ is neither a star or a geometric caterpillar.

By Theorem \ref{urrutia}, the geometric graph $(S^\prime )^c$ contains a plane spanning tree, that is a plane spanning tree $R^\prime$ of $G$ all of whose edges lie in $X^\prime = E\left( G\right) \setminus E\left( S^\prime \right)  =\left( X  \setminus \{ x\} \right) \cup \{ y \} $. Since $X$ is heterochromatic and $c \left( x \right) =  c\left( y \right)$, the set $X^\prime$ and the tree $R^\prime$  are also heterochromatic.

In an analogous way, one may find a heterochromatic plane spanning tree of $G$ in the case where $S$ is a star.

\end{proof}

\begin{theorem}

Let $P$ be a set of  $n\geq3$  points in general position in the plane and let $G$ be the complete geometric graph with vertex set $P$. If $i   \left( P \right)  \leq 1$ and $H = H\left( G\right)$ is the hypergraph with vertex set $E\left( G\right)$ whose hyperedges are the sets of edges of plane spanning trees of $G$, then  $h_{c} \left( H\right)  = \nu \left( H\right) - \tau \left( H\right )+2$.

\end{theorem} 

\begin{proof}
By Theorem \ref{geometricas1},  $h_{c} \left( H\right)   \leq  {n \choose 2}  - n + 2$. For the case where $i  \left( V \left( G \right) \right) = 0$, Lemma \ref{transversal1} asserts that the boundary of $CH\left( P\right)$ is a double transversal of hyperedges of $H$. Hence  $\tau\left(H\right) \leq n$, and by the general lower bound,  $h_{c} \left( H\right)   \geq \nu \left( H\right) - \tau \left( H\right) + 2 = {n \choose 2}  - n + 2$.

If $i  \left( V \left( G \right) \right) = 1$, by Theorem  \ref{geometricas2}, $h_{c} \left( H\right)   \leq  {n \choose 2}  - n + 1$. By  Lemma  \ref{transversal2}, $G$ contains a double transversal of plane spanning trees with $n+1$ edges; therefore $\tau\left(H\right) \leq n+1$. By the general lower bound, $h_{c} \left( H\right)   \geq \nu \left( H\right) - \tau \left( H\right) + 2 = {n \choose 2}  - (n+1) + 2 = {n \choose 2}  - n + 1$.

\end{proof}

\section{Matroids}

For a simple connected graph $G$ with at least three vertices, we denote by $ \gamma\left(G \right)$ the smallest integer $k$ for which there is a set of $k$ edges of $G$ whose removal brakes $G$ into three connected components.

J. Arocha and V. Neumann-Lara  \cite{AN} proved  that if $G$ is a simple connected graph with $m\geq2$ edges and $c$ is a colouring of the edges of $G$ with exactly $m  - \gamma\left(G \right) +2$ colours, then $G$ has an heterochromatic spanning tree.  We generalise this result for matroids.

\begin{theorem}
\label{matroides}

Let $M$ be a matroid with $m$ elements and rank at least 2 and let $\tau$ denote the size of a smallest double transversal of bases of $M$. If $c$ is a colouring of $M$ with exactly $m  - \tau +2$ colours, then $M$ has a heterochromatic basis.

\end{theorem}

\begin{proof}

Denote by $E$ the ground set of $M$ and let $X \subset E$ be a heterochromatic set with $m  - \tau +2$
elements. The complementary set $Y=E\setminus X$ cannot be a double
transversal of bases of $M$ since $\left\vert Y\right\vert =$ $\left\vert
E\setminus X\right\vert = \tau -2$. Hence, there is a
basis $R$ of $M$ that meets $Y$ in at most one element. 

If $R$ is not heterochromatic, then $R$ contains an element $x \in X$ such that 
$c\left( x\right) =c\left( y\right) $, where $y  \neq x$ is the unique element in $
R\cap Y$. 

Let $Z=Y\cup \left\{ x\right\} $, as $\left\vert Z\right\vert =\left\vert
Y\right\vert +1 = \tau -1$, set $Z$ is not a double
transversal of bases of $M$ either. This implies that there is a basis $S$ of $M$
that intersects $Z$ in at most one element.

Assume $S$ is not heterochromatic in which case $\left\vert S\cap
Z\right\vert =1$ and $\left\vert S\setminus Z\right\vert =\left\vert
S\right\vert -1$. Since $\left\vert R\setminus \left\{ x,y\right\}
\right\vert =\left\vert R\right\vert -2=\left\vert S\right\vert -2=\left\vert S\setminus Z\right\vert - 1$, there
is an element $z \in S\setminus Z\subset X$ such that $\left( R\setminus
\left\{ x,y\right\} \right) \cup \left\{ z\right\} $ is an independent set
of $M$. This implies that either $x$ or $y$ must lie in the unique circuit
contained in $R\cup \left\{ z\right\} $ and therefore either $\left( R\cup
\left\{ z\right\} \right) \setminus \left\{ x\right\} $ or $\left( R\cup
\left\{ z\right\} \right) \setminus \left\{ y\right\} $ is a
basis of $M$, a heterochromatic basis of $M$.

\end{proof}

As an immediate consequence we obtain the following:

\begin{corollary}

Let $M$ be a matroid with rank $r\geq 2$.  If $H=H\left( M\right)$ is the hypergraph whose vertices and hyperedges are the elements of the ground set and the bases of $M$, respectively, then $h_{c} \left( H\right)  = \nu \left( H\right) - \tau \left( H\right )+2$.

\end{corollary}

\section{Final remarks and acknowledgements}

We conjecture that if $H$ is the hypergraph associated as in section 2 to any complete geometric graph $G$ with $n\geq3$ vertices in general position, then $h_{c} \left( H\right)  = \nu \left( H\right) - \tau \left( H\right )+2$. Moreover, we conjecture that if  $G$ is a complete geometric graph with $n\geq3$ vertices in general position, then there is a double transversal of plane spanning trees of $G$ with $n+i  \left( V \left( G \right) \right)$ edges and that if $c$ is a colouring of the edges of  $G$ with ${n \choose 2}  - (n+i  \left( V \left( G \right) \right)) + 2$ colours, then $G$ has a heterochromatic plane spanning tree.

We thank Marc Noy for suggesting a way to improve the readability of one of our original proofs.

\end{document}